\documentclass{SCIS2018}
\usepackage{epstopdf}
\usepackage{picins}
\usepackage{diagbox}
\usepackage{graphics} 

\begin{document}

\newtheorem{Remark}{Remark}
\newtheorem{Theorem}{Theorem}
\newtheorem{Assumption}{Assumption}
\newtheorem{Property}{Property}
\newtheorem{Corollary}{Corollary}
\ArticleType{RESEARCH PAPER}
\Year{2018}
\Month{}
\Vol{61}
\No{}
\DOI{}
\ArtNo{}
\ReceiveDate{}
\ReviseDate{}
\AcceptDate{}
\OnlineDate{}

\title{A Parameter Formula Connecting PID and ADRC}{A Parameter Formula Connecting PID and ADRC}

\author[1,2]{Sheng Zhong}{}
\author[1,2]{Yi Huang}{{yhuang@amss.ac.cn.}}
\author[1,2]{Lei Guo}{}

\AuthorMark{Sheng Zhong, et al}

\AuthorCitation{Sheng Zhong, Yi Huang, Lei Guo}


\address[1]{Key Lab. of Systems and Control, Academy of Mathematics and Systems Science,\\
        Chinese Academy of Sciences, Beijing {\rm 100190}, China}
\address[2]{School of Mathematical Sciences, University of Chinese Academy of Sciences, Beijing, {\rm 100049}, China}

\abstract{
This paper presents a parameter formula connecting the well-known proportional-integral-derivative (PID) control and the active disturbance rejection control (ADRC).
On the one hand, this formula gives a quantitative lower bound to the bandwidth of the extended state observer (ESO) used in ADRC, implying that the ESO is not necessarily of high gain.
On the other hand, enlightened by the design of ADRC, a new PID tuning rule is provided, which can guarantee both strong robustness and nice tracking performance of the closed-loop systems under the PID control.
Moreover, it is proved that the ESO can be rewritten as a suitable linear combination of the three terms in PID, which can give a better estimate for the system uncertainty than the single integral term in the PID controller.
Theoretical results are verified also by simulations in the paper.
}


\maketitle

\section{Introduction}
Despite of the remarkable progress of modern control theory over the past sixty years, it is widely recognized that the classical proportional-integral-derivative (PID) control is by far the most widely and successfully used controller in engineering systems\cite{Tariq}. However, it has also been pointed out that most of the practical PID loops are poorly tuned, and there is strong evidence that PID controllers remain poorly understood\cite{Aidan}.
Therefore, as mentioned in\cite{Astrom}, better understanding of the PID control may considerably improve its widespread practice, and so contribute to better product quality. Recently, some theoretical investigations on the global convergence of the PID controller for a basic class of nonlinear uncertain systems are given \cite{Guo1} \cite{zhang}, where some necessary and sufficient conditions for the selection of the PID parameters are provided. These results have rigorously demonstrated in theory that the PID controller does have large-scale robustness with respect to both the uncertain nonlinear structure of the plant and the selection of the controller parameters.

On the other hand, the active disturbance rejection control(ADRC), which was originally proposed by Han in 1998\cite{han}, has attracted more and more attention in both theory and applications\cite{xue} -\cite{gao}.
This is largely because of its unique ideas and superior performance, which were readily translated into something valuable in engineering practice: the ability in dealing with a vast range of uncertainties and great transient response \cite{hy}.
Thus, the high level of robustness and the superior transient performance turn out to be the most valuable characteristics of ADRC to make it an appealing solution in dealing with real world control problems.
However, the research on the theoretical analysis for ADRC was progressing haltingly, especially, on how to tune the ADRC parameters to achieve satisfactory performance of the closed-loop system under practical restrictions.

In this paper, we will provide a new parameter formula for the design of PID controller, which is derived from the inherent but rarely noticed relationship between PID and ADRC.
This formula is found to be beneficial for the design of both PID and ADRC.
On the one hand, this formula gives a quantitative lower bound for the bandwidth of the extended state observer (ESO) used in ADRC, implying that the ESO is not necessarily of high gain, thanks to the parameter manifold provided recently in \cite{Guo1} \cite{zhang} for the selection of PID parameters for nonlinear uncertain systems.
On the other hand, enlightened by the structure of the reduced-order ESO in ADRC, a new and concrete tuning rule for PID parameters is found from the unbounded parameter manifold given in \cite{Guo1} \cite{zhang}, which can guarantee global convergence, strong robustness and nice tracking performances, both for the transient phase and the steady state.
Moreover, we will show that the ESO actually corresponds to a suitable linear combination of the proportional-integral-derivative terms in PID, and will also demonstrate that the ESO can give better estimates for the system uncertainty than the single integral term of PID controller.
Our theoretical results are also verified by some numerical simulations.


The rest of the paper is organized as follows.
The detailed problem description is presented in Section 2.
Section 3 introduces the main results of this paper.
Some simulation verifications of the theoretical analysis are given in Section 4.
Finally, the conclusion is presented in Section 5.

\section{Problem Description}
Consider the following second-order nonlinear uncertain system
\begin{equation}\label{sys}
\begin{cases}
\dot {x}_1= x_2,\\
\dot {x}_2=f(x_1,x_2,t)+u(t),\\
\end{cases}
\end{equation}
where $ (x_1,x_2) \in R^2$ is the system state vector and can be measured,
$u(t)$ is the control input, $f(x_1,x_2,t)\in R$ is an unknown nonlinear function of the state $ (x_1,x_2)$ and time $t$.

The control objective is to make the controlled variable $x_1$ track a given bounded reference signal $y^*(t)$, which satisfies
$$\lim\limits_{t \to \infty }{y^*(t)=y^{**}}, \lim\limits_{t \to \infty }{\dot y^*(t)=0},\lim\limits_{t \to \infty }{\ddot y^*(t)=0},$$
where $\dot y^*(t), \ddot y^*(t)$ are the first and second derivatives of $y^*(t)$, respectively, and $y^{**}$ is a constant.

To have a nice transient control performance, we introduce the following desired transient process to be tracked by $x_1(t)$, which is shaped from $y^*(t)$ by a stable linear filter:
\begin{equation}\label{r}
\ddot r=-2 c_r \dot r-c_r^2 (r-y^*(t)),r(0)=x_1(0),\dot r(0)= x_2(0),
\end{equation}
where $c_r$ is a parameter for tuning the speed of the transient process.

In this paper, the classical PID controller for the system (\ref{sys}) is described as follows:
\begin{equation}\label{upid1}
u_{pid}=-{k_{p}}({x_1}-r)-{k_{d}}({x_2}-\dot r)-k_{i}\int_{0}^{t} ({x_1(\tau)}-r(\tau)) d\tau+\ddot r ,
\end{equation}
where $k_p,k_d,k_i$ are the controller parameters to be discussed in the paper.

On the other hand, according to the idea of ADRC, $f$ can be viewed as the total disturbance of the system and treated as an extended state of the system to be estimated by an extended state observer (ESO) so that it can be compensated for in time.

Since the state $x_2$ is measurable, the following reduced-order ESO can be designed \cite{hy}:
 \begin{equation}\label{eso}
\begin{cases}
\dot \xi=-\omega_o \xi -\omega_o^2 x_2-  \omega_o u,\quad \xi (0) = -\omega_o x_2 (0),\\
 \hat f=\xi+\omega_o x_2, \\
\end{cases}
\end{equation}
where $\hat f$ is the estimation of the total disturbance $f(x_1,x_2,t)$, $\hat f(0)=0,$ and $\omega_o$ is the parameter of ESO to be tuned.

Then, the corresponding ADRC law for tracking the transient process $r(t)$ can be designed as \cite{hy}
\begin{equation}\label{u}
u=-{k_{ap}}({x_1}-r)-{k_{ad}}({x_2}-\dot r)-\hat f+\ddot r,\quad k_{ap}>0,k_{ad}>0,
\end{equation}
where $k_{ap},k_{ad}$ are two controller parameters to be tuned.
In the ADRC law (\ref{u}), 
the term $-\hat f$, which is an estimate of $f$, tries to compensate for the total disturbance, and $\ddot r$ is a feedforward term.
Thus, the ADRC (\ref{u}) can be regarded as an adaptive pole-placement control with given closed-loop poles determined by $k_{ap}$ and $k_{ad}.$

Substituting the equation (\ref{u}) into (\ref{eso}) gives
\begin{equation}\label{a1}
\hat f={\omega_o  k_{ad}}({x_1}-r)+ \omega_o ({x_2}-\dot r)+\omega_o  k_{ap} \int_{0}^{t} ({x_1(\tau)}-r(\tau)) d\tau .
\end{equation}
Therefore, the ADRC law (\ref{u}) can be rewritten as£º
\begin{equation}\label{u2}
u=-({k_{ap}}+\omega_o k_{ad}) ({x_1}-r)-({k_{ad}}+\omega_o)({x_2}-\dot r)-\omega_o k_{ap}\int_{0}^{t} ({x_1(\tau)}-r(\tau)) d\tau+\ddot r .
\end{equation}
Comparing (\ref{u2}) with (\ref{upid1}) shows that if set
\begin{equation}\label{a2}
k_p = k_{ap}+\omega_o k_{ad}, \quad k_d=k_{ad}+\omega_o,\quad k_i = \omega_o k_{ap},
\end{equation}
then the PID (\ref{upid1}) is the same as the ADRC law (\ref{u}).

The above simple parameter formula (\ref{a2}), which connects PID and ADRC, is quite meaningful. It suggests that

1. The main results provided in \cite{Guo1} and \cite{zhang} on the selection of PID parameters for guaranteeing the global asymptotic stability of the closed-loop system, may be used to find quantitative lower bounds for the parameters $(k_{ap},k_{ad},\omega_o)$ of ADRC (\ref{eso}) and (\ref{u}).
In the next section, this quantitative lower bound for ADRC will be firstly given.
This result will show that the parameters of ADRC are not necessary of high gain.
Moreover, theoretical analysis will demonstrate that the performance of the closed-loop system may be improved by tuning the parameter $\omega_o.$

2. The formula (\ref{a2}) provides a new and concrete tuning rule for PID parameters rather than taken arbitrarily from a given unbounded parameter manifold as in \cite{Guo1} \cite{zhang}.
Furthermore, when the parameters $(k_{p},k_{d},k_i)$ of PID are tuned by the formula (\ref{a2}), the suitable linear combination of the P part $\omega_o k_{ad} (x_1-r)$, the D part $\omega_o  (x_2-\dot r)$ and the I part $\omega_o k_{ap}\int_{0}^{t} ({x_1(\tau)}-r(\tau)) d\tau$ has a good function of estimating the unknown $f$.
In the classical PID controller, the integral part has certain capability to estimate $f$ at least in the steady state.
In the next section, it will be proved that the output $\hat f$ of ESO (\ref{eso}), which is the combination of the three terms of P-I-D (\ref{a1}), has a better capability for estimating the dynamic process of an unknown function $f$.


\section{Main Results}
Before presenting the main results, we introduce a definition for a class of unknown nonlinear functions $f$.
Define the following function space:
\begin{equation}\label{flim}\begin{split}
 \mathcal{F}=&\Big \{ f\in C^1 (R^2 \times R^+) \Big | f(x_1,x_2,t)=h(x_1,x_2)+w(t), \Big | \frac{\partial h}{\partial x_1}\Big |\leq L_1,
  \Big | \frac{\partial h}{\partial x_2}\Big |\leq L_2, \Big |w(t) \Big |\leq L_3,\\
  & \Big |\dot w (t)\Big |\leq L_3, \lim\limits_{t \to \infty }{w(t)}  \text{ exists},  \forall x_1,x_2 \in R, \forall t \in R^+ \Big \},
\end{split}\end{equation}
where $L_1,L_2,L_3$ are positive constants, and $C^1 (R^2 \times R^+)$ denotes the space of all functions from $R^2 \times R^+$ to $R$ which are
continuous in $(x_1,x_2)$ and $t$,
with continuous partial derivatives with respect to $(x_1,x_2)$.


\subsection{A lower bound to the gain of the ESO with guaranteed performance}
As suggested by the formula (\ref{a2}) and the manifold provided in \cite{Guo1} \cite{zhang} for the selection of PID parameters, a quantitative lower bound to the parameter $\omega_o$ of the ESO (\ref{eso}) may be obtained.
To this end, we denote

$$ \Omega=\big \{\omega \in R \big | n_0 \omega^4 + n_1  \omega^3 + n_2 \omega^2 + n_3 \omega + n_4 = 0 \big \},$$
where $n_0 = k_{ad}^2$ and
\begin{equation}\begin{split}
&  n_1=2 k_{ad} [k_{ad} (k_{ad}-L_2)-L_1], \\
&n_2 = 2 k_{ad} (k_{ap}-L_1)(k_{ad}-L_2)+ [k_{ad} (k_{ad}-L_2)-L_1]^2-L_2^2k_{ap},\\
&n_3=2 [k_{ad} (k_{ad}-L_2)-L_1](k_{ap}-L_1)(k_{ad}-L_2)-L_2^2k_{ap}(k_{ad}-L_2), \\
&n_4=(k_{ap}-L_1)^2(k_{ad}-L_2)^2.
\end{split}\end{equation}
We can now obtain a lower bound $\omega_o^*$ to the ESO parameter $\omega_o$ as follows.
\begin{equation}\label{omegalim}\begin{split}
&{\bar \omega_o=
\begin{cases}
0,& \Omega = \varnothing  \quad or  \quad \max \{\Omega \} \leq 0  ,\\
\max \{\Omega \},& \max \{\Omega \} > 0,
\end{cases}}\\
&\omega_o^*=\max \begin{Bmatrix}0,\frac{L_1- k_{ap}}{k_{ad}},L_2-k_{ad},\bar \omega_o\end{Bmatrix},\\
\end{split}\end{equation}
where $\varnothing$ represents the empty set.
\begin{Remark}
From (\ref{omegalim}), it can be seen that $\omega_o^*$ only depends on the constants $L_1,L_2,k_{ap},k_{ad}$ and is irrelevant to the disturbance $w(t$), initial values and the reference signal $y^*(t).$
\end{Remark}

The following theorem shows that $\omega_o^*$ can indeed serve as a lower bound to the ESO parameter $\omega_o.$

\begin{Theorem}\label{Thm1}
\emph{\textbf{(Tracking Performance)}}.
Consider the ADRC controlled nonlinear uncertain system (\ref{sys}),(\ref{eso}) and (\ref{u}), where the nonlinear unknown function $f \in  \mathcal{F}.$
Then, for any given $L_1,L_2,k_{ap},k_{ad}$, the closed-loop tracking error will satisfy
$$\lim\limits_{t \to \infty }{x_1(t)=y^{**}}, \lim\limits_{t \to \infty }{x_2(t)=0},$$
for any initial value $(x_1(0),x_2(0)) \in R^2 $ and any $y^{**}$, as long as the ESO parameter $\omega_o > \omega_o^*$.
\end{Theorem}

The proof of Theorem \ref{Thm1} is given in Appendix A.

Theorem \ref{Thm1} gives a tuning method of ADRC which makes the closed-loop system achieve global asymptotic stability ultimately.
It can be seen from Theorem \ref{Thm1} that the lower bound to the parameter of the ESO (\ref{eso}), i.e., $\omega_o^*$, can be calculated through (\ref{omegalim}).
This result indicates that the parameter of ESO is not necessarily of high gain.

The next Theorem will further show that the tracking performance may be improved by tuning the parameter $\omega_o > \omega_o^*$.

Denote the tracking error as $e(t)= r(t) - x_1(t)$ and the estimation error as $e_f=\hat f-f(x_1,x_2,t)$.

\begin{Theorem}\label{Thm2}
\emph{\textbf{(Pole-placement Performance)}}.
Consider the ADRC controlled nonlinear uncertain system (\ref{sys}),(\ref{eso}) and (\ref{u}), where the nonlinear unknown function $f \in  \mathcal{F}.$
Then, there exist positives constants $\eta_1,\eta_2 ,$ which depend on
$(e_f(0),k_{ap},$ $k_{ad},L_1,L_2,L_3,\dot r, \ddot r)$, such that for all $\omega_o > \omega_o^*$, the closed-loop equation of the pole-placement ADRC has the following property:
\begin{equation}\label{thm21}
|\ddot e(t) + k_{ad} \dot e(t) + k_{ap}e(t)| =|e_f(t)| \leq \eta_1 e ^{-\omega_o t} + \frac{\eta_2}{\omega_o},\quad t \geq 0 .
\end{equation}
\end{Theorem}

The proof of Theorem \ref{Thm2} is given in Appendix A.

Note that $\ddot e(t) + k_{ad} \dot e(t) + k_{ap}e(t)=0$ is the ideal performance for the tracking error of pole-placement control.
Therefore, Theorem \ref{Thm2} describes the distance between the real closed-loop performance and the ideal one.
From (\ref{thm21}), the dynamic response can be divided into two parts.
The first part $\eta_1 e ^{-\omega_o t}$, which is related to the initial value of the estimation error $e_f$, can be rapidly tuned to zero by the parameter $\omega_o$.
The second part $\frac{\eta_2}{\omega_o}$ decreases with $\omega_o$.
Thus, by increasing the parameter $\omega_o$, the dynamic response of the closed-loop system can be made close to the ideal trajectory.

\subsection{ A new tuning rule for PID}
In this section, a new tuning rule for PID controller is proposed, which can guarantee the robustness as well as nice tracking performance of the closed-loop system.

According to the parameter formula (\ref{a2}) and Theorem \ref{Thm2}, a new tuning rule for the PID law (\ref{upid1}) is given as follows.
\begin{equation}\begin{split}\label{tuning}
&k_p = k_{ap}+\omega_o k_{ad}, \quad k_d=k_{ad}+\omega_o,\quad k_i = \omega_o k_{ap},\\
& k_{ap}>0,\quad k_{ad}>0,\quad \omega_o > \omega_{o}^*.
\end{split}
\end{equation}

Under the tuning rule (\ref{tuning}), the PID controller (\ref{upid1}) is equivalent to the ADRC (\ref{eso}), (\ref{u}).
That is to say, the properties of the closed-loop system (\ref{sys}), (\ref{upid1}) are the same as those of the closed-loop system (\ref{sys}),(\ref{eso}) and (\ref{u}).
Thus, according to the advantages of ADRC, the PID controller defined by (\ref{upid1}) and (\ref{tuning}) also has the ability to timely estimate and compensate for the disturbances and uncertainties, so that the closed-loop system has guaranteed strong robustness and superior tracking performance.

The following corollary can be directly obtained based on Theorem \ref{Thm2}.

\begin{Corollary}\label{coro}
Consider the PID controlled nonlinear uncertain system (\ref{sys}),(\ref{upid1}) and (\ref{tuning}), where the nonlinear unknown function $f \in  \mathcal{F}.$
Then, there exist positives constants $\eta_1,\eta_2 ,$ which are the same as in Theorem~\ref{Thm2}, such that for any given $L_1,L_2,k_{ap},k_{ad}>0$, any initial value $(x_1(0),x_2(0)) \in R^2 $ and any setpoint $y^{**}$,the closed-loop system has the following properties whenever $\omega_o > \omega_o^*$:
\begin{itemize}
\itemindent 2.8em
\item[(1)] $\lim\limits_{t \to \infty }{x_1(t)=y^{**}}, \lim\limits_{t \to \infty }{x_2(t)=0}$.
\item[(2)] $|\ddot e(t) + k_{ad} \dot e(t) + k_{ap}e(t)| \leq \eta_1 e ^{-\omega_o t} + \frac{\eta_2}{\omega_o},\quad t \geq 0.$
\end{itemize}
\end{Corollary}


To further elaborate on the nice performance of the ESO, we note that the integral term $\hat f_I=k_i \int_{0}^{t} ({x_1(\tau)}-r(\tau)) d\tau$ of PID controller (\ref{upid1}) is usually regarded to have the ability to eliminate the constant disturbance, while in the ADRC frame, the ESO (\ref{eso}) is known to have the ability to timely estimate the dynamic disturbance.
The following theorem compares the estimation property of the integral term of PID controller (\ref{upid1}) with that of the ESO (\ref{eso}) in the frequency domain.

Define $ e_{f_I}(t)=\hat f_I -f(x_1,x_2,t)$.
Denote $E_f(s), E_{f_I} (s),\hat F_I (s)$ and $\hat F (s)$ as the Laplace transforms of $e_f(t), e_{fI}(t),\hat f_I(t)$ and $\hat f(t)$, respectively.
It can be obtained from Theorem \ref{Thm1} that the unknown $f$ on the system trajectories is bounded, thus the Laplace transform of $f$ exists.
Denote $F(s)$ as the Laplace transform of the unknown $f$, and let $G_{e_f}(s),G_{e_{f_I}}(s)$ be the transfer functions from $F(s)$ to $E_f(s)$ and $E_{f_I} (s)$, respectively.

\begin{Theorem}\label{pro1}
Consider the system (\ref{sys}), (\ref{eso}) and (\ref{u}) and the system (\ref{sys}), (\ref{upid1}), which are connected by the formula (\ref{tuning}).
For any $f \in  \mathcal{F},$ the integral term of PID controller (\ref{upid1}) and the ESO (\ref{eso}) have the following properties when $\omega_o > \omega_o^*$:
\begin{itemize}
\item[(1)]For any $\omega$, ${ \frac{|G_{e_f} (i\omega)|}{|G_{e_{f_I}} (i\omega)|}  }<1.$ Moreover, $\lim\limits_{t \to \infty }{\frac{e_f(t)}{e_{f_I}(t)}} =  \frac{k_{ap} }{k_{ap}+\omega_o k_{ad} }.$
\item[(2)]$\hat F_I (s)=\frac{ k_{ap}}{s^2+k_{ad} s + k_{ap}}\hat F (s) .$
\end{itemize}
\end{Theorem}

The proof of Theorem \ref{pro1} is given in Appendix A.

From Theorem \ref{pro1}(1), it can be seen that the steady estimation error of the ESO (\ref{eso}) is smaller than that of the integral term for the total disturbances $f$ at any frequency. 
Moreover, it can be seen that
the ratio of the steady estimation error is $\frac{k_{ap} }{k_{ap}+\omega_o k_{ad} }$, which decreases with the increase of either $\omega_o$ or $k_{ad}$.
The result (2) of Theorem \ref{pro1} shows that phase-lag of the response of the ESO (\ref{eso}) is smaller than that of the integral term, particularly for rapidly varying disturbances.

\begin{Remark}
From the formula (\ref{a1}), (\ref{tuning}) and Theorem \ref{pro1}, it can be concluded that the P-term and the D-term of PID controller also contribute to the estimation and compensation for the disturbances, rather than the single I-term.
This seems to be a somewhat striking property that has not been revealed in the investigation of the classical PID before.
\end{Remark}

\section{Simulations}
In this section, some simulations are presented to verify the main results of the paper.

In the simulations, the unknown nonlinear $f$ can be one of following cases:
\begin{equation*}\begin{split}
&C1: f(t)=2x_1+6sinx_2+1,\qquad C2: f(t) =5cos x_1 + 2x_2 -2,\\
&C3: f(t)=3 x_1 +2 x_2 -w_1 (t), \qquad C4: f(t)=6 cosx_1 sin x_2 -w_2(t),\\
\end{split}\end{equation*}
where
$${
w_1(t)=
\begin{cases}
sin(t),  & \text{if $t < 4s$, } \\
sin(4), & \text{else,}
\end{cases}
\quad w_2(t)=
\begin{cases}
cos(t),  & \text{if $t < 4s$, } \\
cos(4), & \text{else.}
\end{cases}}$$
The initial values of the state are $x_1(0)=0,x_2(0)=0,$ and the reference signal is $y^*(t)=2.$
Then, the desired transient process $r(t)$ is designed as follows:
\begin{equation}\label{r1}
\ddot r=-2 c_r \dot r-c_r^2 (r-2),c_r=5 ,r(0)=0,\dot r(0)= 0.
\end{equation}

According to Theorem \ref{Thm1}, it can be calculated that $\omega_o^*=6$.
In the simulations, the parameters in the formula (\ref{tuning}) are chosen as $k_{ap}=4,k_{ad}=4,\omega_o=10.$

\begin{figure}[!t]
\centering
\begin{minipage}[c]{0.48\textwidth}
\centering
\includegraphics [scale=0.5] {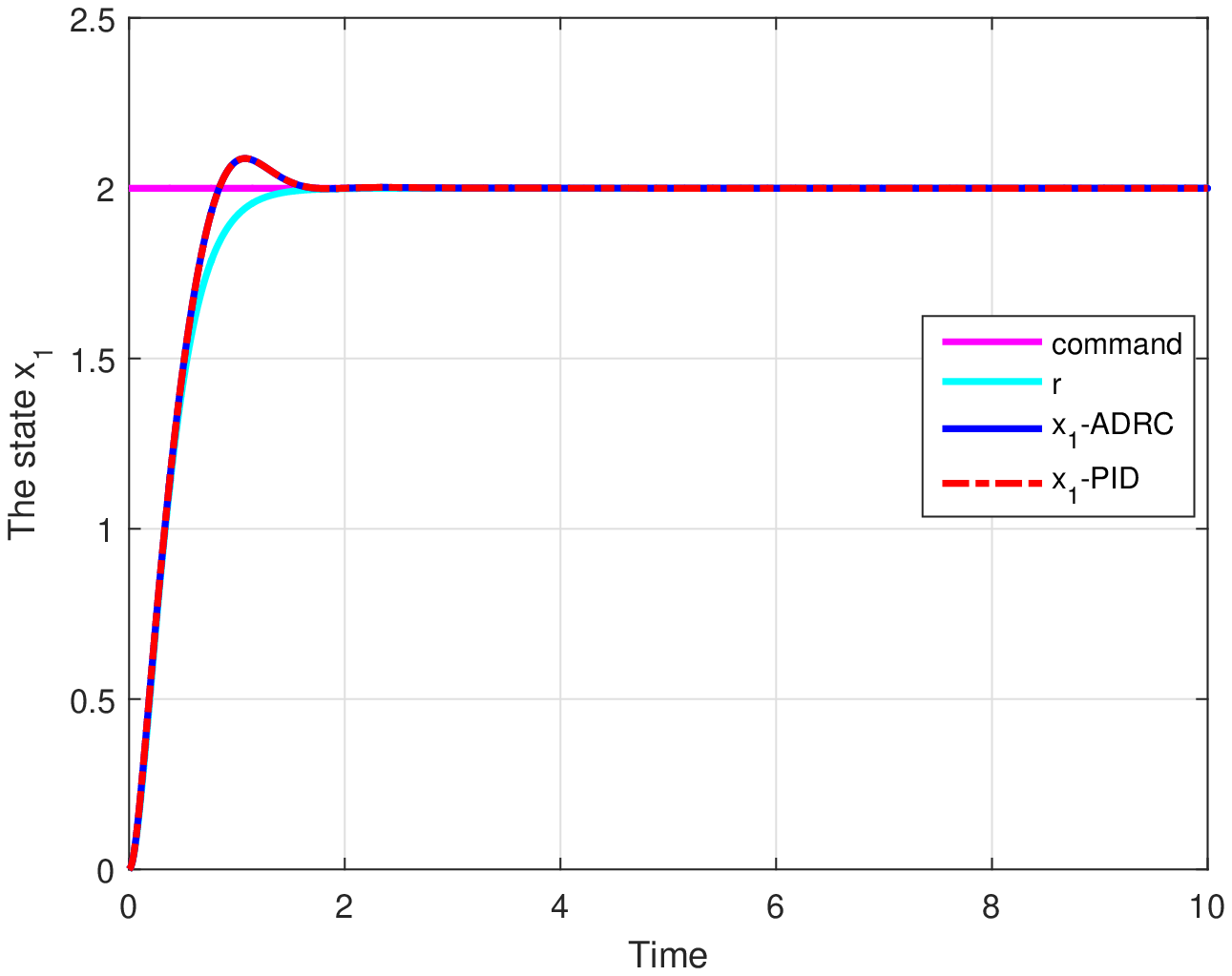}
\caption*{(a)}
\end{minipage}
\hspace{0.02\textwidth}
\begin{minipage}[c]{0.48\textwidth}
\centering
\includegraphics [scale=0.5] {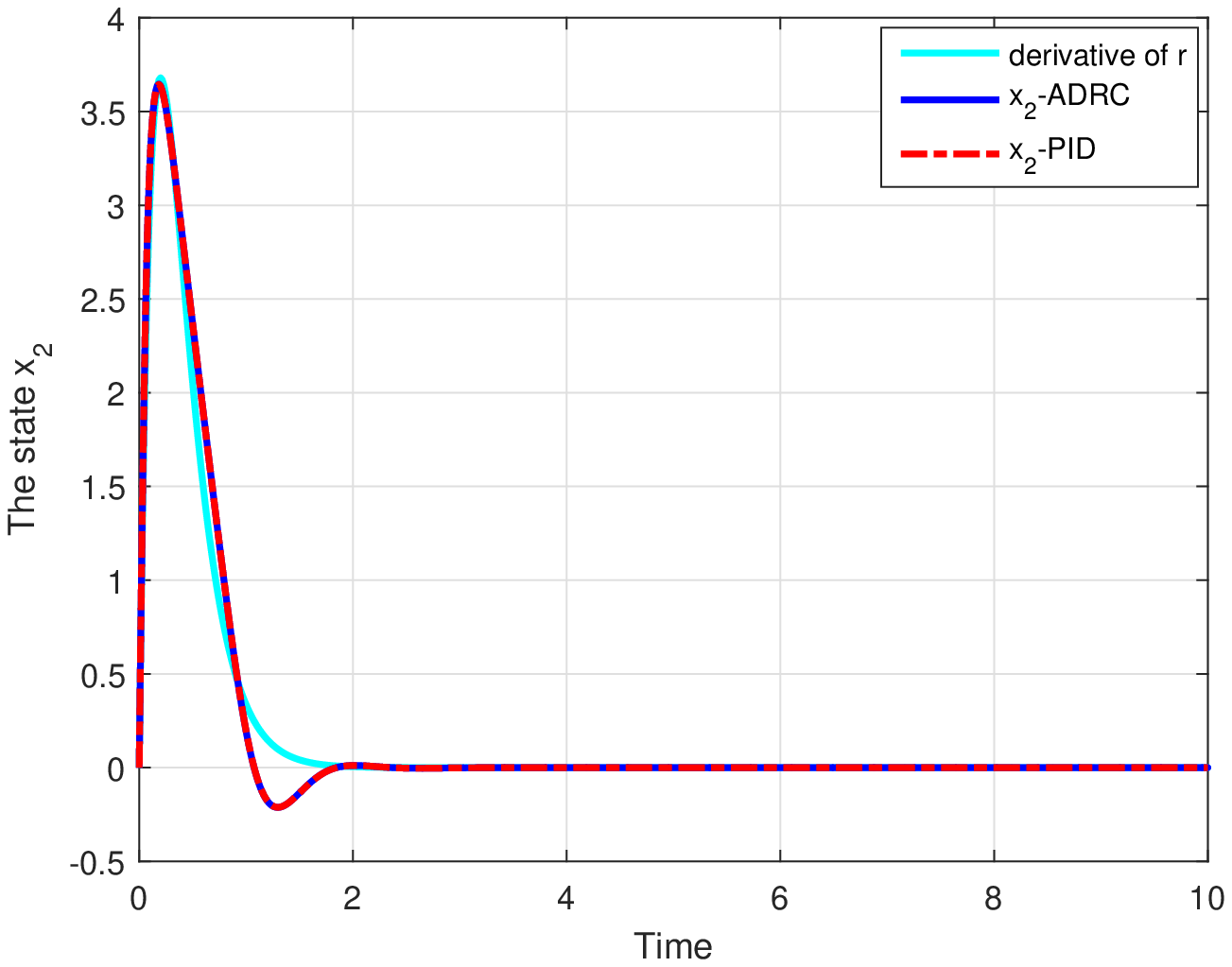}
\caption*{(b)}
\end{minipage}
\caption{(a) The response curves of the state $x_1$; (b) The response curves of the state $x_2$.}
\label{fig1}
\end{figure}

\begin{figure}[!t]
\centering
\begin{minipage}[c]{0.48\textwidth}
\centering
\includegraphics [scale=0.5] {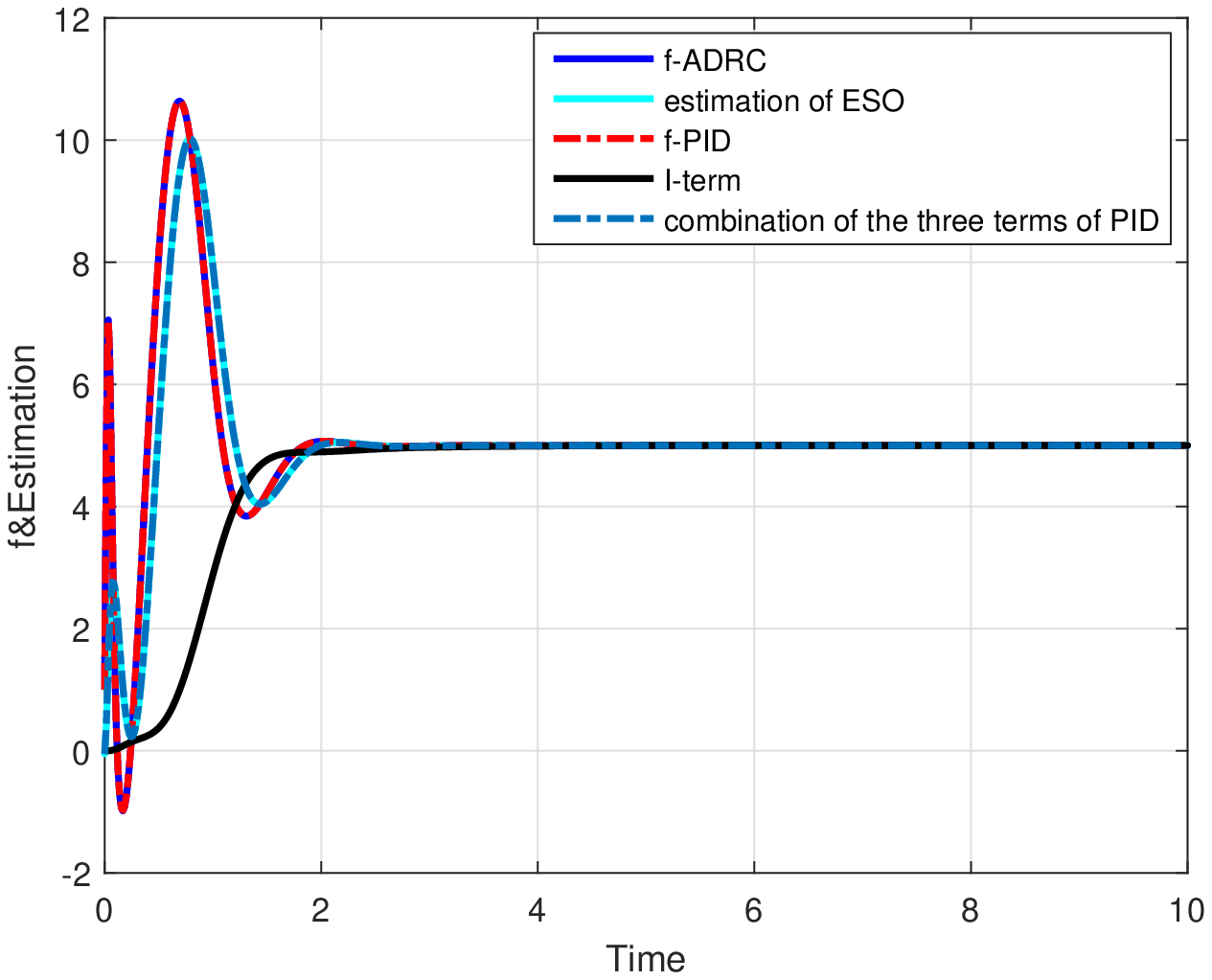}
\caption*{(a)}
\end{minipage}
\hspace{0.02\textwidth}
\begin{minipage}[c]{0.48\textwidth}
\centering
\includegraphics [scale=0.5] {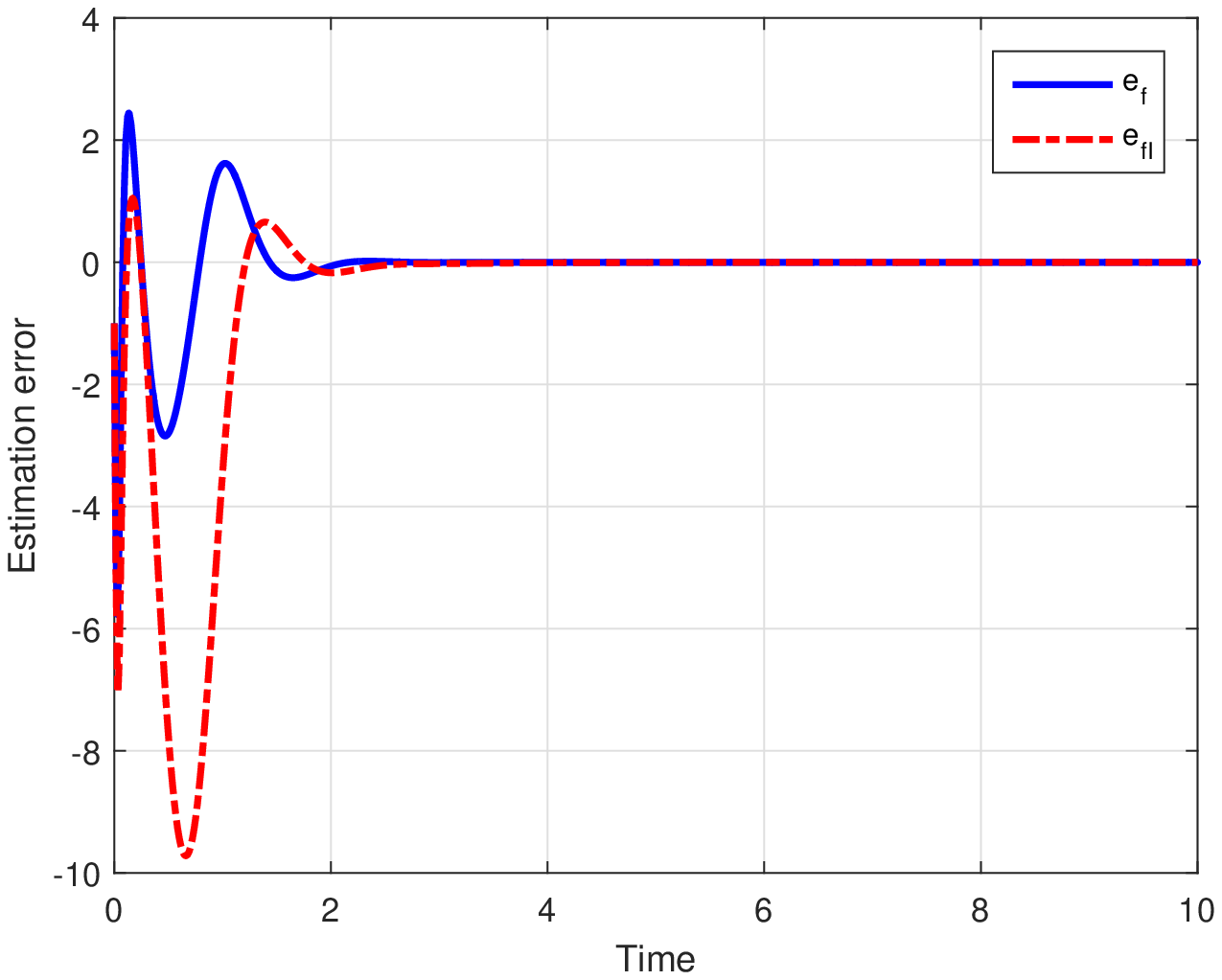}
\caption*{(b)}
\end{minipage}
\caption{(a) The estimations of the disturbance $f$; (b) The estimation errors of the disturbance $f$.}
\label{fig2}
\end{figure}

The simulation results in the $C1$ case are shown in Figures \ref{fig1} and \ref{fig2}.

Figure \ref{fig1} is the response curves of the state $(x_1,x_2)$ based on the ADRC (\ref{eso}) and (\ref{u}) (the blue line), and the PID controller (\ref{upid1}) and (\ref{tuning}) (the red dash line).
It is shown that under the parameter formula (\ref{tuning}), the closed-loop system (\ref{sys}) and (\ref{upid1}) and the closed-loop system (\ref{sys}),(\ref{eso}) and (\ref{u}) have the same dynamic responses.

Figure \ref{fig2}(a) are the estimations of the unknown $f$ based respectively on the ESO (\ref{eso}), the I-term of PID controller (\ref{upid1}) and (\ref{tuning}), and the combination of the three terms of P-I-D of PID controller (\ref{a1}).
In Figure \ref{fig2}(a), the blue line represents $f$ of the closed-loop system (\ref{sys}),(\ref{eso}) and (\ref{u}), and the red dash line represents $f$ of the closed-loop system (\ref{sys}), (\ref{upid1}) and (\ref{tuning}).
Figure \ref{fig2}(b) are the estimation errors of $f$ based on the ESO (\ref{eso}) (the blue line) and the I-term of PID controller (\ref{upid1}) and (\ref{tuning}) (the red dash line), respectively.

Figure \ref{fig2}(a) shows that compared to the I-term, the ESO (\ref{eso}) can track the unknown disturbance more quickly.
Moreover, it also verifies that the combination of the three terms of P-I-D (\ref{a1}) has the same capability for estimating the unknown $f$ as that of the ESO (\ref{eso}).
From Figure \ref{fig2}(b), it can be seen that the estimation error of the ESO (\ref{eso}) is smaller than that of the I-term, although both gradually tend to zero.

\begin{figure}[!t]
\centering
\begin{minipage}[c]{0.48\textwidth}
\centering
\includegraphics[scale=0.5]{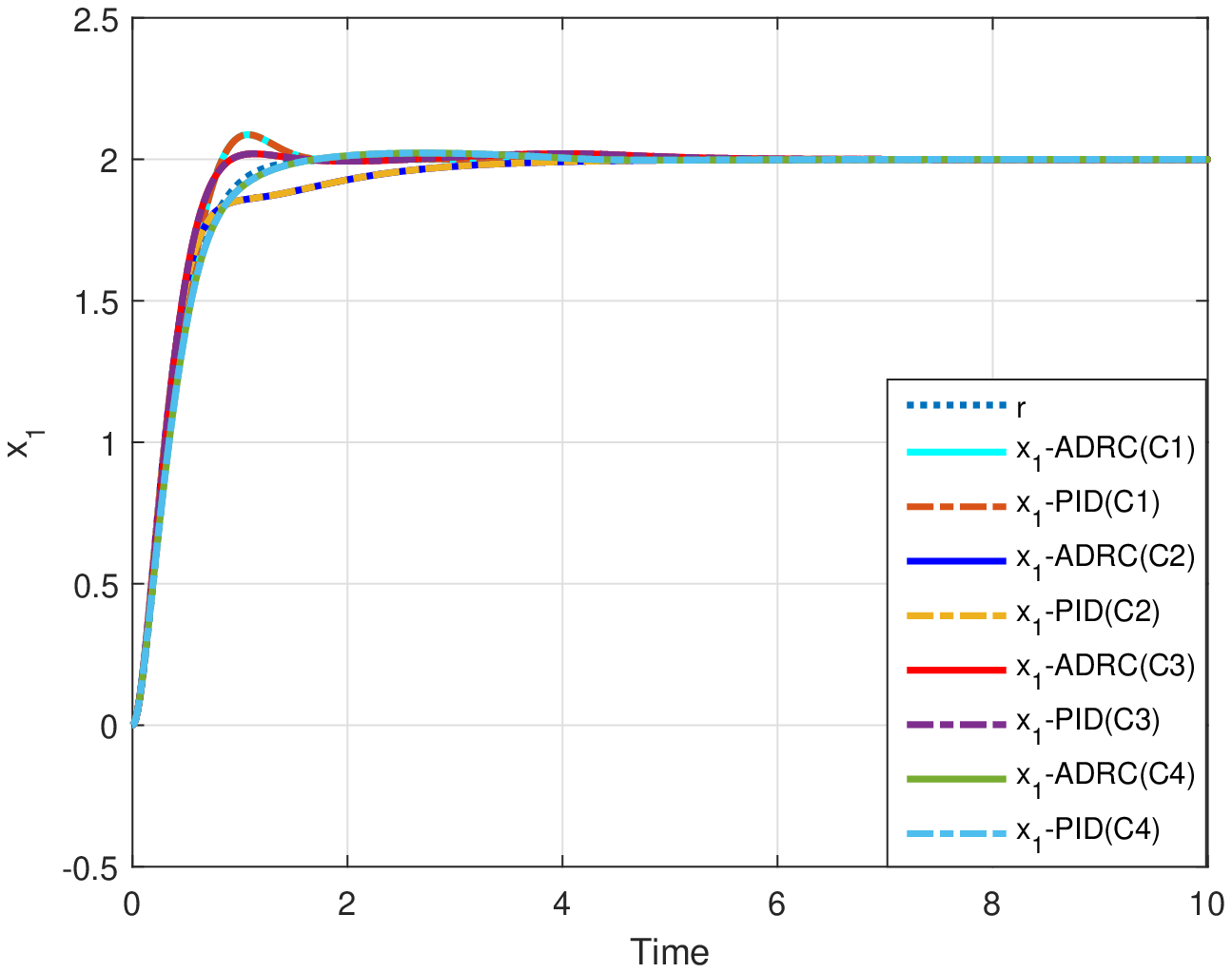}
\end{minipage}
\hspace{0.02\textwidth}
\begin{minipage}[c]{0.48\textwidth}
\centering
\includegraphics[scale=0.5]{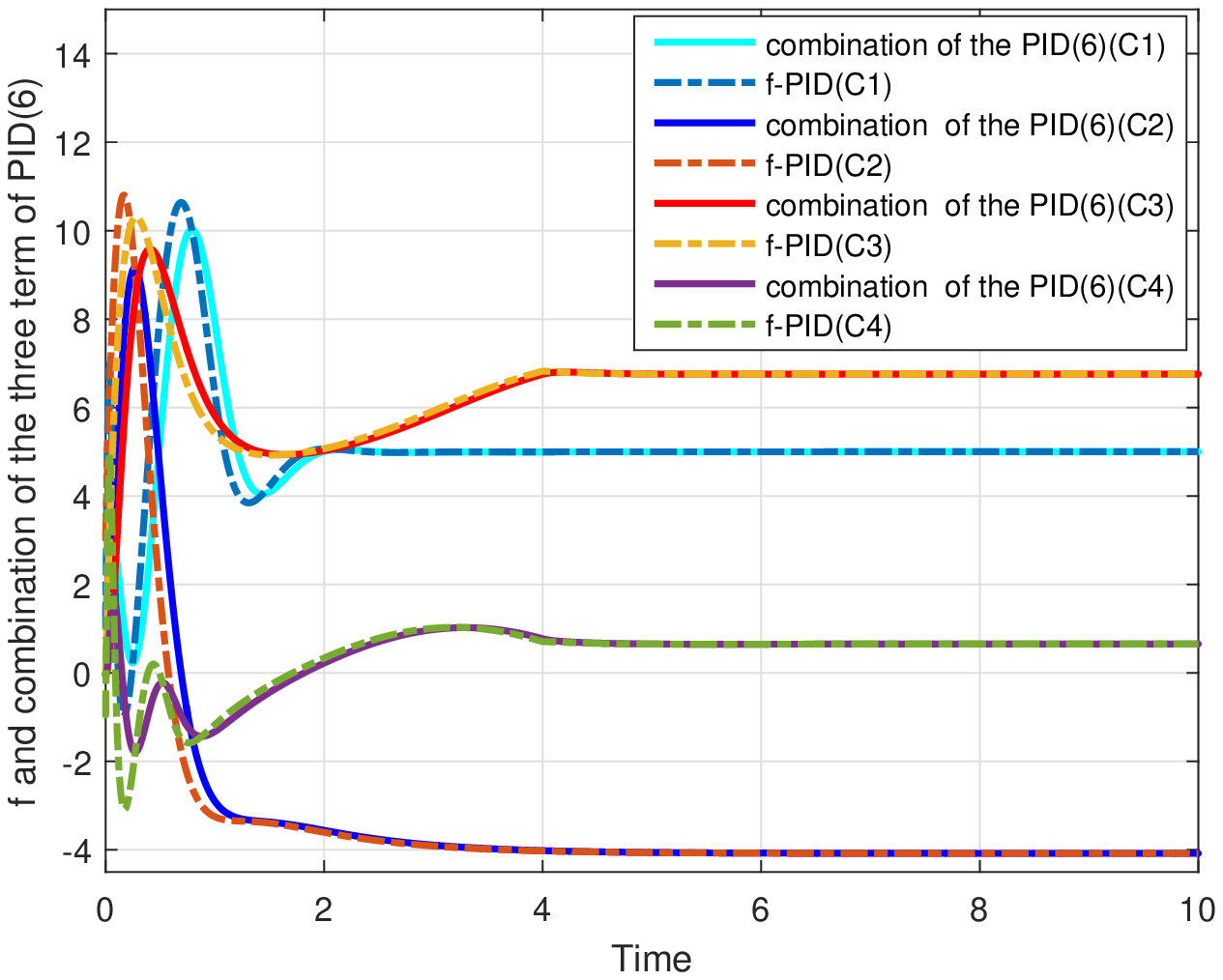}
\end{minipage}\\[3mm]
\begin{minipage}[t]{0.48\textwidth}
\centering
\caption{The response curves of the state $x_1$ under $C1 \sim C4$.}
\label{fig3}
\end{minipage}
\hspace{0.02\textwidth}
\begin{minipage}[t]{0.48\textwidth}
\centering
\caption{The curves of disturbances $f$ and combination of the three terms of P-I-D, under $C1 \sim C4$.}
\label{fig4}
\end{minipage}
\end{figure}

Figure \ref{fig3} are the response curves of the state $x_1$ , based on the ADRC (\ref{eso}) and (\ref{u}), and the PID controller (\ref{upid1}) and (\ref{tuning}) under difference situations $C1 \sim C4$, respectively.
It can be seen that both the PID (\ref{upid1}) and ADRC (\ref{eso}) and (\ref{u}), which are connected by the formula (\ref{tuning}), can deal with a vast range of uncertainties in the sense that the closed-loop system has strong robustness and great tracking performance.



Figure \ref{fig4} are the curves of the unknown $f$ and the combination of the three terms of P-I-D (\ref{a1}) in the cases $C1 \sim C4$, respectively.
It indicates that the combination of the three terms of P-I-D (\ref{a1}), which is equal to the estimation given by the ESO (\ref{eso}), has the ability to timely estimate a large range of the unknown dynamic function $f$.
This is the reason why both the PID (\ref{upid1}) and ADRC (\ref{eso}) - (\ref{u}), tuned according to the formula (\ref{tuning}), have the capability to keep the tracking performance close to the ideal one $r(t)$.


\begin{figure}[!t]
\centering
\begin{minipage}[c]{0.48\textwidth}
\centering
\includegraphics [scale=0.5] {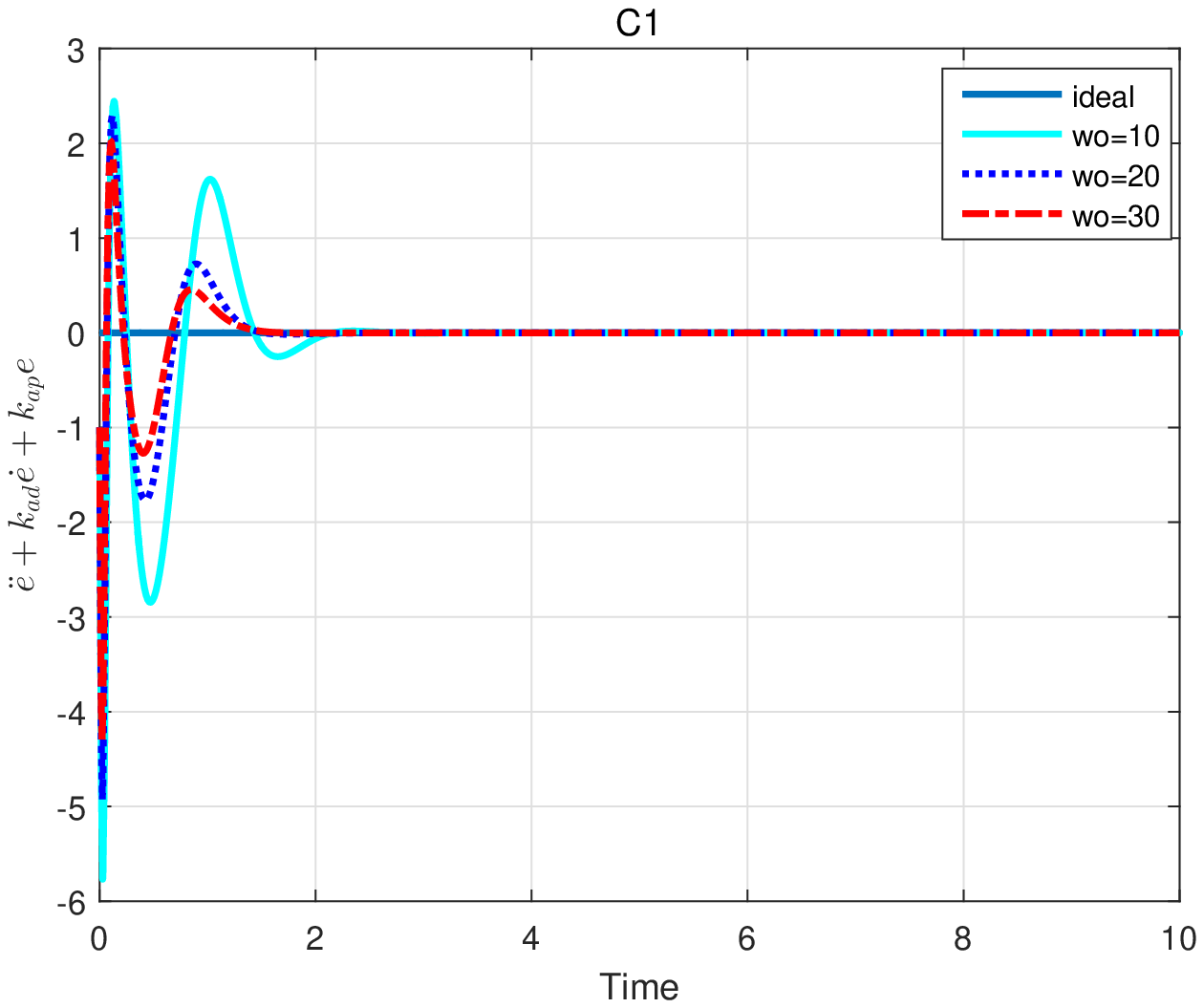}
\caption*{(a)}
\end{minipage}
\hspace{0.02\textwidth}
\begin{minipage}[c]{0.48\textwidth}
\centering
\includegraphics [scale=0.5] {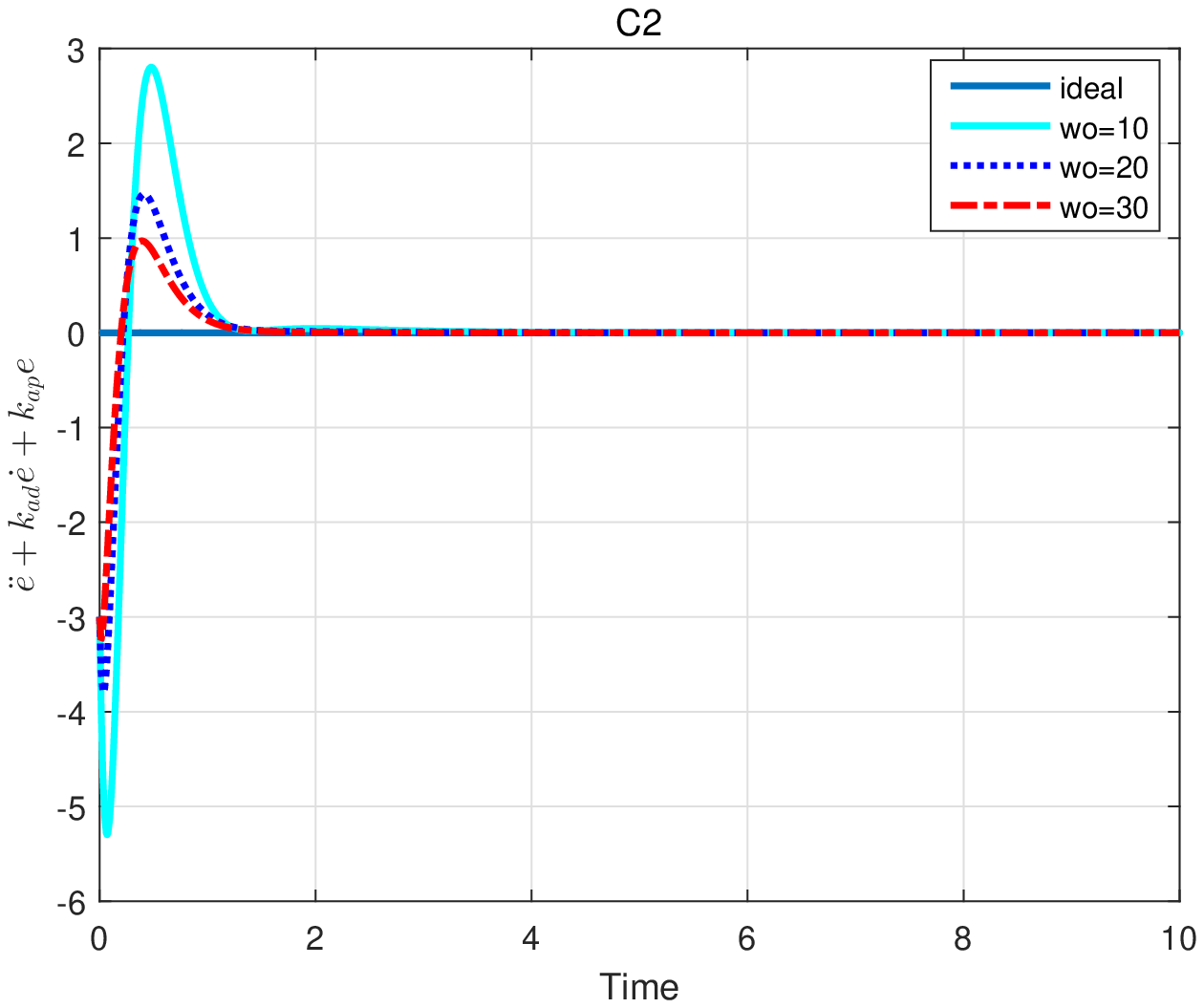}
\caption*{(b)}
\end{minipage}

\begin{minipage}[c]{0.48\textwidth}
\centering
\includegraphics [scale=0.5] {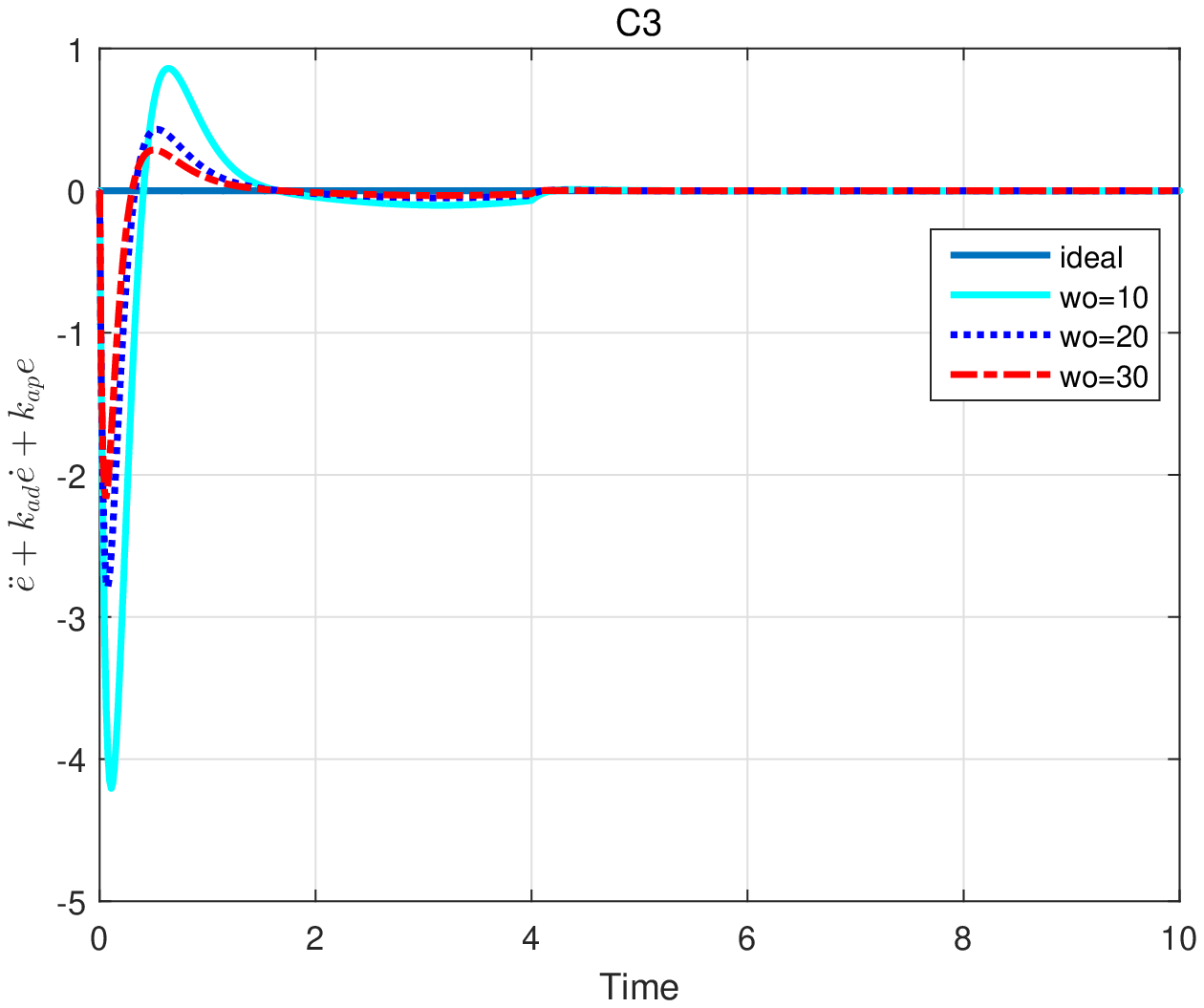}
\caption*{(c)}
\end{minipage}
\hspace{0.02\textwidth}
\begin{minipage}[c]{0.48\textwidth}
\centering
\includegraphics [scale=0.5] {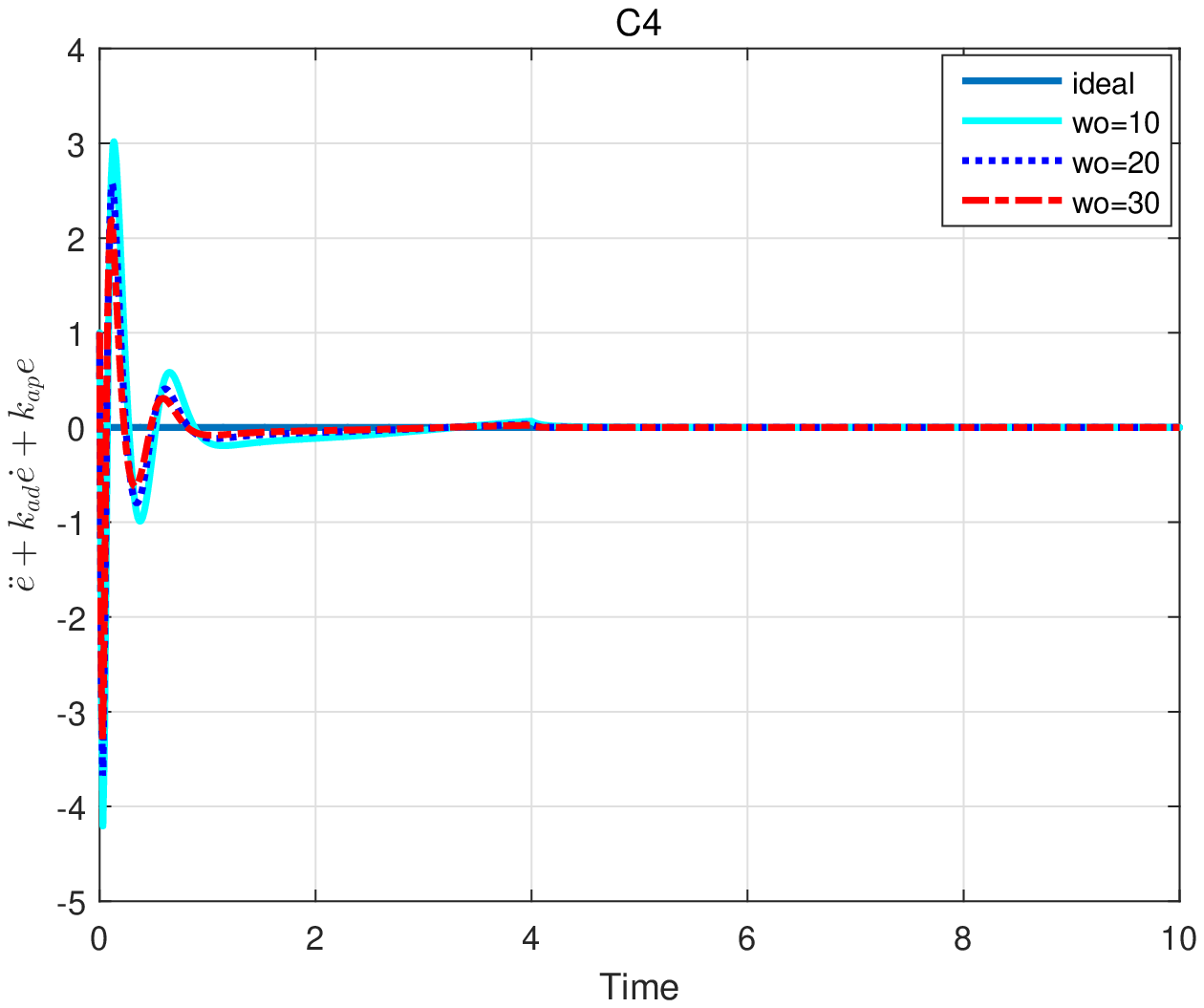}
\caption*{(d)}
\end{minipage}
\caption{(a) The curves of $\ddot e(t) + k_{ad} \dot e(t) + k_{ap}e(t)$ under $C1$; (b) The curves of $\ddot e(t) + k_{ad} \dot e(t) + k_{ap}e(t)$ under $C2$;
(c) The curves of $\ddot e(t) + k_{ad} \dot e(t) + k_{ap}e(t)$ under $C3$; (d) The curves of $\ddot e(t) + k_{ad} \dot e(t) + k_{ap}e(t)$ under $C4$.}
\label{fig5}
\end{figure}

To verify the results of Corollary \ref{coro}, Figure \ref{fig5} is the curves of $\ddot e(t) + k_{ad} \dot e(t) + k_{ap}e(t)$ based on the PID controller (\ref{upid1}), (\ref{tuning}), when the parameter $\omega_o$ varies.
It shows that the increase of the parameter $\omega_o$ will lead to the decrease of $|\ddot e(t) + k_{ad} \dot e(t) + k_{ap}e(t)|$, i.e., the dynamic response of the tracking error $e$ can be tuned close to the ideal trajectory by only increasing $\omega_o$.

\section{Conclusion}
In this paper, a new and simple parameter formula connecting PID and ADRC is discovered and is shown to be able to improve the design of both controllers significantly.
Firstly, a quantitative lower bound $\omega_o^*$ to $\omega_o$, the parameter of the ESO (\ref{eso}), is provided for guaranteeing the global asymptotic stabilizability of the ADRC.
This result shows that the design parameters of the ADRC are not necessary of high gain.
It is further proved that the upper bound for the tracking performance of the closed-loop system can be improved by increasing $\omega_o>\omega_o^*$.
Then, a novel PID controller tuning rule, suggested by the design of ADRC, is provided.
Since the PID controller is equivalent to the ADRC by this tuning rule, the robustness and excellent tracking performance of the closed-loop system can be guaranteed.
Finally, it is demonstrated that the steady estimation error of the ESO (\ref{eso}) is lesser and phase-lag of the response of the ESO (\ref{eso}) is smaller than that of the single integral term of PID controller (\ref{upid1}).
We believe that the tuning formula provided in this paper has wide applicability in practical control systems.
\Acknowledgements{This work was supported by the National Center for Mathematics
and Interdisciplinary Sciences, Chinese Academy of Sciences and the National Natural Science Foundation of China (Grant Nos. 11688101).}



\begin{appendix}
\section{}
{\bf {Proof of Theorem \ref{Thm1}}}

Substituting equation (\ref{u}) into (\ref{eso}), by Laplace transform, we obtain
\begin{equation}\label{lz12}
\hat F(s)={\omega_o  k_{ad}}({X_1(s)}-R(s))+ \omega_o s({X_1(s)}-R(s))+ \frac{\omega_o  k_{ap}({X_1(s)}-R(s))}{s} ,
\end{equation}
where $X_1(s),R(s)$ are the Laplace transforms of the state $x_1(t)$ and the transient process $r(t)$, respectively.
Take the inverse Laplace transform for (\ref{lz12}), we have
\begin{equation}\label{z12}
\hat f={\omega_o  k_{ad}}({x_1}-r)+ \omega_o ({x_2}-\dot r)+\omega_o  k_{ap} \int_{0}^{t} ({x_1(\tau)}-r(\tau)) d\tau .
\end{equation}
Hence, the control law (\ref{u}) can be rewritten as
\begin{equation}\label{upid}
u=-{k_{p}}({x_1}-r)-{k_{d}}({x_2}-\dot r)-k_{i}\int_{0}^{t} ({x_1(\tau)}-r(\tau)) d\tau+\ddot r ,
\end{equation}
where $ k_p = k_{ap}+\omega_o k_{ad}, k_d=k_{ad}+\omega_o, k_i = \omega_o k_{ap}.$

Since $ \lim\limits_{t \to \infty }{w(t)}  \text{ exists}$, which can be denoted by a constant $c$.
Denote $$ e_i(t)= \int_{0}^{t} e(\tau) d\tau + \frac{h(y^{**},0)+c}{k_i}, \quad e_d(t)= \dot e(t),\quad
g(e,e_d)= -h(y^{**}-e,-e_d)+h(y^{**},0).$$
Based on the definition of $\mathcal{F}$, it can be seen that $g\in \mathcal{F}$ and $g(0,0)=0$.
Then the closed-loop system (\ref{sys}) and (\ref{upid}) turns into
\begin{equation}\label{closesys}
\begin{cases}
\dot {e_i}= e,\\
\dot e=e_d,\\
\dot e_d= -k_i e_i -k_p e -k_d e_d +g(e,e_d)+\Delta(t) ,
\end{cases}
\end{equation}
where $\Delta(t)=g(e+y^{**}-r,e_d-\dot r) -g(e,e_d)+c-w(t).$
By the mean value theorem, it can be obtained that, for any $t\in R^+,$ there is
\begin{equation}\label{delta}
g(e+y^{**}-r,e_d-\dot r) -g(e,e_d)=\frac{\partial g}{\partial e} \Big | _ {(\bar e,\bar e_d)}(y^{**}-r)-\frac{\partial g}{\partial e_d} \Big | _ {(\bar e,\bar e_d)} \dot r ,
\end{equation}
where $\bar e = e+\theta (y^{**}-r), \bar e_d = e_d-\theta \dot r, \theta \in (0,1)$.
Since $f\in \mathcal{F}$, it can be deduced that $|\Delta| \leq L_1 |y^{**}-r|+L_2 |\dot r|+c+L_3.$
Moreover, $(0,0,0)$ is an equilibrium of (\ref{closesys}), when $t$ approaches infinity.

Following the analysis in \cite{Guo1} and \cite{zhang}, we denote
$$b(e)=
\begin{cases}
\frac{g(e,0)}{e},& e\neq 0,\\
\frac{\partial g}{\partial e}(0,0),& e=0,
\end{cases}
\quad and \quad
a(e,e_d)=
\begin{cases}
\frac{g(e,e_d)-g(e,0)}{e_d},& e_d\neq 0,\\
\frac{\partial g}{\partial e_d}(e,0),& e_d=0,
\end{cases}
$$
then $g(e,e_d)$ can be expressed as $$g(e,e_d)= b(e) e + a(e,e_d) e_d.$$
By the mean value theorem again and the definition of $\mathcal{F}$, we have $| b(e)| \leq L_1, | a(e,e_d)| \leq L_2$ for all $e,e_d$.

Hence, the closed-loop system (\ref{closesys}) can be rewritten as
\begin{equation}\label{closesys1}
\begin{cases}
\dot {e_i}= e,\\
\dot e=e_d,\\
\dot e_d= -k_i e_i -\phi (e) e - \psi (e,e_d) e_d  +\Delta(t),
\end{cases}
\end{equation}
where $\phi (e)=k_p - b(e), \psi (e,e_d)=k_d - a(e,e_d).$
By the fact that $\omega_o > \frac{L_1- k_{ap}}{k_{ad}}, \omega_o > L_2-k_{ad},$ there exist $\phi (e)\geq k_p-L_1>0$ and $\psi (e,e_d) \geq k_d-L_2>0$.

To construct a Lyapunov function, we consider the following matrix $P$:
\begin{equation}\label{lyapp}
P=\frac{1}{2}\begin{bmatrix} \mu k_i & k_i & \delta\\ k_i & k_p-L_1 + \mu k_d & \mu\\ \delta & \mu & 1 \end{bmatrix}, \quad
\mu = \frac{2((k_p-L_1)(k_d-L_2) + k_i)}{4(k_p-L_1) + L_2^2}
\end{equation}
where $\delta$ satisfies $0<\delta < \frac{2\mu k_i}{(k_p-L_1) + \mu k_d}$.
We will first show that the matrix $P$ is positive definite.

Based on the definition of $\omega_o^*$ and the assumption $\omega_o>\omega_o^*,$ it can be obtained that
\begin{equation}\label{lyapp1}\begin{split}
(k_p-L_1)(k_d-L_2)-k_i>L_2\sqrt{k_i(k_d-L_2)},
\end{split}\end{equation}
thus,
\begin{equation}\label{lyapp2}\begin{split}
&(k_p-L_1)(k_d-L_2)-k_i>0,\\
&[(k_p-L_1)(k_d-L_2)-k_i]^2>L_2^2{k_i(k_d-L_2)}.
\end{split}\end{equation}
Since
\begin{equation}\label{lyapp3}\begin{split}
\mu-k_d+L_2= \frac{-2(k_p-L_1)(k_d-L_2) + 2k_i-L_2^2 (k_d-L_2)}{4(k_p-L_1) + L_2^2}<0,
\end{split}\end{equation}
and
\begin{equation}\label{lyapp4}\begin{split}
4(-\mu+k_d-L_2)(-k_i+\mu(k_p-L_1))-\mu^2 L_2^2= \frac{4[[(k_p-L_1)(k_d-L_2)-k_i]^2-L_2^2{k_i(k_d-L_2)}]}{4(k_p-L_1) + L_2^2}>0,
\end{split}\end{equation}
we have
\begin{equation}\label{lyapp5}\begin{split}
-k_i+\mu(k_p-L_1)>0.
\end{split}\end{equation}
Then, based on (\ref{lyapp3})-(\ref{lyapp5}), the following three inequalities can be verified.
\begin{equation}\label{lyapp6}\begin{split}
\mu k_i>0,
\left|\begin{array}{cccc}
    \mu k_i &    k_i  \\
    k_i &    k_p-L_1 + \mu k_d
\end{array}\right|
=k_i[\mu( k_p-L_1 + \mu k_d )-k_i]>0,
\end{split}\end{equation}
and
\begin{equation}\label{lyapp7}\begin{split}
\left|\begin{array}{cccc}
   \mu k_i & k_i & \delta\\ k_i & k_p-L_1 + \mu k_d & \mu\\ \delta & \mu & 1
\end{array}\right|
>k_i(\mu(k_p-L_1)+\mu^2 k_d -k_i -\mu^3)>0.
\end{split}\end{equation}
Thus, the matrix $P$ is positive definite.

We are now in a position to consider the following Lyapunov function (\cite{Guo1}\cite{zhang}):
\begin{equation}\label{lyap}
V(e_i,e,e_d)=[e_i,e,e_d]P[e_i,e,e_d]^T + \int_{0}^{e} (L_1-b(s))s d s.
\end{equation}
Since $0\leq \int_{0}^{e} (L_1-b(s))s d s \leq L_1 e^2$, from (\ref{lyap}), we have
\begin{equation}\label{PF2_7}
 [e_i,e,e_d]P[e_i,e,e_d]^T\leq V(e_i,e,e_d)\leq [e_i,e,e_d]P_0[e_i,e,e_d]^T,
 \end{equation}
where
$$
P_0=\frac{1}{2}\begin{bmatrix} \mu k_i & k_i & \delta\\ k_i & k_p+L_1 + \mu k_d & \mu\\ \delta & \mu & 1 \end{bmatrix}.
$$
It can be deduced that
\begin{equation}\label{PF2_8}
\lambda_{\min}(P) \left\Vert [e_i,e,e_d] \right\Vert^2 \leq {V} \leq\lambda_{\max}(P_0) \left\Vert {[e_i,e,e_d]}\right\Vert^2,
\end{equation}
where $\lambda_{\min}(\cdot)$ and $\lambda_{\max}(\cdot)$ are the minimum and maximum eigenvalues of the corresponding matrix, respectively.

The time derivative of $V(e_i,e,e_d)$ along the trajectories of (\ref{closesys1}) is
\begin{equation}\label{dlyap}
\dot V(e_i,e,e_d)=-[e_i,e,e_d]W(e,e_d) [e_i,e,e_d]^T+ [\delta , \mu, 1][e_i,e,e_d]^T\Delta ,
\end{equation}
where 
$$W(e,e_d)=\begin{bmatrix} \delta k_i & \frac{\delta \phi (e)}{2} & \frac{\delta \psi (e,e_d)}{2}\\ \frac{\delta \phi (e)}{2} & -k_i + \mu \phi (e) & -\frac{\mu a(e,e_d)+\delta}{2}\\ \frac{\delta \psi (e,e_d)}{2} & -\frac{\mu a(e,e_d)+\delta}{2} &
-\mu+\psi(e,e_d) \end{bmatrix}.$$

Denote
$$Q(e,e_d)=\begin{bmatrix} -k_i + \mu \phi (e) & -\frac{\mu a(e,e_d)}{2}\\ -\frac{\mu a(e,e_d)}{2} &
-\mu+\psi(e,e_d) \end{bmatrix}.$$
From (\ref{lyapp3})-(\ref{lyapp5}), it is easy to verify that $Q(e,e_d)$ is positive definite. Since $\phi (e), \psi(e,e_d)$ are bounded, based on the same analysis as in \cite{zhang}, there exists a constant $\delta^*$, such that the matrix $W(e,e_d)$ is positive definite when $\delta<\delta^*$. Thus $\delta$ can be defined by $\delta< \min\{\frac{2\mu k_i}{k_p-L_1 + \mu k_d},\delta^*\}.$

From (\ref{dlyap}), we have
\begin{equation}\label{dlyap1}
\dot V\leq -\lambda_{\min}(W)\|[e_i,e,e_d]\| ^2+ c_0 \|[e_i,e,e_d]\| |\Delta |\leq -c_1 V +c_2 \sqrt V| \Delta| ,
\end{equation}
where $c_0 = \max\{\delta , \mu, 1\}, c_1 = \frac{\lambda_{\min}(W)}{\lambda_{\max}(P_0)}, c_2 = \frac{c_0}{\sqrt{\lambda_{\min}(P)}}.$
Because $|\Delta|$ is bounded, there exists a constant $M_0$, such that $|\Delta| \leq M_0$.
Then, it can be obtained that
\begin{equation}\label{dlyap2}
\sqrt V \leq  e ^{\frac{-c_1 t}{2}}\sqrt {V(e_i(0),e(0),e_d(0))} + \frac{c_2 M_0}{c_1} (1-e ^{\frac{-c_1 t}{2}})\leq M_1,
\end{equation}
where $M_1$ is a constant.
Since $ \lim\limits_{t \to \infty }{r(t)=y^{**}}, \lim\limits_{t \to \infty }{\dot r(t)=0}$ and $ \lim\limits_{t \to \infty }{w(t)=c},$
we know that $ \lim\limits_{t \to \infty }{\Delta =0}$. Thus, for any $\varepsilon>0$, there exists $T>0$, such that for any $t>T$, there is
$\dot V\leq -c_1 V + \varepsilon$.

In conclusion, $ \lim\limits_{t \to \infty }{x_1(t)=y^{**}}, \lim\limits_{t \to \infty }{x_2(t)=0}$.
The proof of Theorem \ref{Thm1} is complete.
\hfill $\blacksquare$

$$\quad $$
{\bf {Proof of Theorem \ref{Thm2}}}

Since the dynamic equations of $e$ and $e_d$ can be written as follows:
\begin{equation}\label{closesys2}
\begin{cases}
\dot e=e_d,\\
\dot e_d= -k_{ap} e -k_{ad} e_d  +e_f,
\end{cases}
\end{equation}
we get $|\ddot e(t) + k_{ad} \dot e(t) + k_{ap}e(t)|= |e_f(t)|$.
Denote
$$
P_2=\begin{bmatrix} \frac{1+k_{ap}}{2k_{ad}} + \frac{k_{ad}}{2k_{ap}}& \frac{1}{2k_{ap}} \\ \frac{1}{2k_{ap}} & \frac{1+k_{ap}}{2k_{ad}k_{ap}}\end{bmatrix}.
$$
Consider the following Lyapunov function
\begin{equation}\label{lyape}
V_2(e,e_d)=\begin{bmatrix}e&e_d\end{bmatrix} P_2 \begin{bmatrix} e\\e_d \end{bmatrix}.
\end{equation}
The time derivative of $V_2 (e,e_d)$ along the trajectories of (\ref{closesys2}) is
\begin{equation}\label{lyape1}
\dot V_2=-e^2 -e_d^2 + \begin{bmatrix}e&e_d\end{bmatrix} \begin{bmatrix}  \frac{1}{k_{ap}} \\ \frac{1+k_{ap}}{k_{ad}k_{ap}}\end{bmatrix}e_f
\leq  - \frac{V_2}{\lambda_{\max}(P_2)}+\frac {\max(\frac{1}{k_{ap}},\frac{1+k_{ap}}{k_{ad}k_{ap}})\sqrt V_2|e_f|} {\sqrt{\lambda_{\min}(P_2)}}.
\end{equation}
From (\ref{lyape1}), it can be seen that
$$\sqrt {V_2} \leq \frac {\max(\frac{1}{k_{ap}},\frac{1+k_{ap}}{k_{ad}k_{ap}})\lambda_{\max}(P_2)\mathop{sup}|e_f|} {\sqrt{\lambda_{\min}(P_2)}} ,$$ i.e., $$\| [e,e_d] \| \leq \frac {\max(\frac{1}{k_{ap}},\frac{1+k_{ap}}{k_{ad}k_{ap}})\lambda_{\max}(P_2)\mathop{sup}|e_f|} {{\lambda_{\min}(P_2)}}.$$
Since when $\sqrt {V_2} > \frac {\max(\frac{1}{k_{ap}},\frac{1+k_{ap}}{k_{ad}k_{ap}})\lambda_{\max}(P_2)\mathop{sup}|e_f|} {\sqrt{\lambda_{\min}(P_2)}} ,$ we know that $\dot V_2<0.$

From the equation (\ref{eso}), the estimation error $e_f$ satisfies the following equation
\begin{equation}\label{ez}
\dot e_f = -\omega_o e_f - \dot f,
\end{equation}
where $$\dot f= \frac{\partial f}{\partial x_2} k_{ap} e +(\frac{\partial f}{\partial x_2} k_{ad}-\frac{\partial f}{\partial x_1})e_d-\frac{\partial f}{\partial x_1}e_f+\frac{\partial f}{\partial x_1 }\dot r+\frac{\partial f}{\partial x_2} \ddot r +\frac{\partial f}{\partial t}.$$
Since $f\in \mathcal{F}$, we have
\begin{equation}
|\dot f| \leq \gamma_1 |e| + \gamma_2 |e_d| +L_2 |e_f|+ \gamma_3,
\end{equation}
where $\gamma_1=L_2 k_{ap},\gamma_2=|L_2 k_{ad}-L_1|,\gamma_3=L_1\dot r+L_2 \ddot r+L_3.$

Consider the following Lyapunov function:
\begin{equation}\label{lyapez}
V_1(e_f)=\frac{1}{2}e_f^2,
\end{equation}
the time derivative of $V_1 (e_f)$ along the trajectories of (\ref{ez}) is
\begin{equation}\label{lyapez1}
\dot V_1=-\omega_o e_f^2- e_f \dot f \leq -\omega_o e_f^2 +|e_f| |\dot f|\leq -(\omega_o-L_2) e_f^2 +|e_f|(\gamma_1 |e| + \gamma_2 |e_d| + \gamma_3).
\end{equation}

Denote $$\gamma_4=\frac {\max(\frac{1}{k_{ap}},\frac{1+k_{ap}}{k_{ad}k_{ap}})\lambda_{\max}(P_2)} {{\lambda_{\min}(P_2)}},\quad  \omega_{o1}^*=L_2+\gamma_4(\gamma_1+\gamma_2)+1.$$
Next, it will be proved that when $\omega_o > \omega_{o1}^*, \mathop{sup}|e_f| \leq max \{ e_f(0),\gamma_3\}.$
From (\ref{lyapez1}), it can be seen that
\begin{equation}\label{lyapez11}
\dot V_1< -(\gamma_4(\gamma_1+\gamma_2)+1) e_f^2 +|e_f|(\gamma_4(\gamma_1+\gamma_2)\mathop{sup}|e_f| + \gamma_3).
\end{equation}
When $|e_f| > max \{ e_f(0),\gamma_3\},$ there is $\dot V_1 <0,$ thus, $|e_f| \leq max \{ e_f(0),\gamma_3\}.$
Moreover, there exists a constant $\gamma_5$, which does not depend on $\omega_o$, such that, $| \dot f| \leq \gamma_5.$

If $\omega_o^* \leq \omega_o \leq \omega_{o1}^*,$ then by Theorem~\ref{Thm1} and the equation (\ref{dlyap1}), we know that there exists a constant $\gamma_6$, such that $\gamma_6=\mathop{sup} \limits_{\omega_o} | \dot f|$.

Denote $M_{ \dot f}=max\{ \gamma_5,\gamma_6\}.$
Based on the above analysis, it can be deduced that
\begin{equation}\label{lyapez2}
\sqrt V_1 \leq  e ^{-\omega_o t}\sqrt {V_1(e_f(0))} + \frac{\sqrt 2 M_{ \dot f}}{2\omega_o} (1-e ^{-\omega_o t}),
\end{equation}
thus
\begin{equation}\label{lyapez3}
|e_f| \leq  e ^{-\omega_o t}|e_f(0)| + \frac{ M_{ \dot f}}{\omega_o} (1-e ^{-\omega_o t}) \leq \eta_1 e ^{-\omega_o t} + \frac{\eta_2}{\omega_o} ,
\end{equation}
where $\eta_1= |e_f(0)| ,\eta_2 = M_{ \dot f},$ which are irrelevant to $\omega_o$.
Thus, (\ref{thm21}) is obtained and the proof of Theorem \ref{Thm2} is complete.
\hfill $\blacksquare$


$$\quad $$
{\bf {Proof of Theorem \ref{pro1}}}

Based on the proofs of Theorems \ref{Thm1} and \ref{Thm2}, when $\omega_o > \omega_o^*$, both $f$ and $\dot f$ are bounded.
Since the PID controller defined by (\ref{upid1}) and (\ref{tuning}) is equivalent to the ADRC (\ref{eso}) and (\ref{u}), the total disturbance $f$ of the closed-loop system defined by (\ref{sys}), (\ref{upid1}) and (\ref{tuning}) is the same as that of the closed-loop system (\ref{sys}), (\ref{eso}) and (\ref{u}).
Based on the equation (\ref{closesys2}), we have
\begin{equation}\label{epez1}
\ddot e +k_{ad} \dot e +k_{ap} e = e_f.
\end{equation}
Take the Laplace transform for the equation (\ref{epez1}), we get
\begin{equation}\label{epez2}
E(s)= \frac{1}{s^2+k_{ad}s+k_{ap}} E_f (s),
\end{equation}
where $E(s)$ is the Laplace transform of $e(t)$ and $ E_f (s)$ is the Laplace transform of $e_f(t)$.
Based on the equation (\ref{ez}), we know that
\begin{equation}\label{epez21}
\dot e_{f_I}=-k_i e + \dot e_f +\omega_o e_f.
\end{equation}
Take the Laplace transform for the equation (\ref{epez21}), it can be obtained from (\ref{epez2}) that
\begin{equation}\label{epez22}
{E_{f_I} (s)} = \frac{s^2 + k_d s + k_p }{s^2+k_{ad} s + k_{ap} } {E_f (s)} ,
\end{equation}
where $ E_{f_I} (s)$ is the Laplace transform of $e_{f_I}(t)$.
Take the Laplace transform for the equation (\ref{ez}), it follows that
\begin{equation}\label{epez3}
E_f (s) =G_{e_f}(s) F(s),G_{e_f}(s)=\frac{s}{s+\omega_o },
\end{equation}
where $ F(s)$ is the Laplace transform of $f.$
Thus,
\begin{equation}\label{epez222}
E_{f_I} (s) = G_{e_{f_I}}(s)F(s),G_{e_{f_I}}(s)=\frac{s^3 + k_d s^2 + k_p s}{(s+\omega_o)(s^2+k_{ad} s + k_{ap})}.
\end{equation}

From the equation (\ref{epez3}) and (\ref{epez222}), we obtain
\begin{equation}\label{epez4}
{ \frac{|G_{e_f} (i\omega)|^2}{|G_{e_{f_I}} (i\omega)|^2}  }=\frac{(k_{ap}- \omega^2)^2 + k_{ad}^2 \omega^2 }{(k_{p}- \omega^2)^2 + k_{d}^2\omega^2 }<1.
\end{equation}
Since $\lim\limits_{t \to \infty } {\frac {e_f(t)}{e_{f_I}(t)}}=\lim\limits_{s \to 0 }{\frac{s E_f(s)}{sE_{f_I}(s)}}$, it is easy to see that
$$\lim\limits_{t \to \infty }{\frac{e_f(t)}{e_{f_I}(t)}} =  \frac{k_{ap} }{k_{ap}+\omega_o k_{ad} }.$$
Therefore, Theorem \ref{pro1} (1) is obtained.

Based on the equations (\ref{sys}) and (\ref{u}), it can be obtained that
\begin{equation}\label{epez7}
\dddot {\hat f_I}= -k_d {\ddot {\hat {f_I}}} -k_p {\dot {\hat {{f}_{I}}}} -k_i \hat f_I + k_i  f.
\end{equation}
Take the Laplace transform for the equation (\ref{epez7}), we have
\begin{equation}\label{epez8}
\hat F_I (s) =\frac{\omega_o k_{ap}}{(s+\omega_o)(s^2+k_{ad} s + k_{ap})}  F(s),
\end{equation}
where $ \hat F_I (s)$ is the Laplace transform of $\hat f_I(t)$ and $ F(s)$ is the Laplace transform of $ f.$

Based on the equation (\ref{eso}), the dynamical equation of $\hat f$ can be written as follows
\begin{equation}\label{z121}
\dot {\hat {f}} = -\omega_o (\hat f -  f).
\end{equation}
Take the Laplace transform for the equation (\ref{z121}), we have
\begin{equation}\label{epez9}
\hat F(s) =\frac{\omega_o}{s+\omega_o}F(s),
\end{equation}
where $ \hat F(s)$ is the Laplace transform of $\hat f(t)$.
Thus,
$$\hat F_{I}(s)=\frac{ k_{ap}}{s^2+k_{ad} s + k_{ap}}\hat F(s),  $$
which is Theorem \ref{pro1} (2).
Hence, the proof of Theorem \ref{pro1} is complete.
\hfill $\blacksquare$

\end{appendix}

\end{document}